
\voffset=-0.4in

\input amstex
\magnification =\magstep 1
\documentstyle{amsppt}

\NoBlackBoxes
\hoffset=4mm

\def\H{{\Cal H}}

\def\N{{\Bbb N}}

\topmatter

\title An Operator Hilbert space without the Operator Approximation Property
\endtitle

\rightheadtext{A Hilbertian space without the OAP}

\author Alvaro Arias
\endauthor

\address Alvaro Arias
\newline Division of Mathematics and Statistics, 
The University of Texas at San Antonio, San
Antonio, TX 78249, U.S.A.\endaddress
\email arias\@math.utsa.edu\endemail

\abstract 
We use a technique of Szankowski [S] to
construct an operator Hilbert space that does not
have the operator approximation property  
\endabstract

\subjclass Primary 46B28
Secondary 46B20, 47D15
\endsubjclass

\endtopmatter

\document
\baselineskip18pt
\head Introduction and Preliminaries \endhead

A Banach space $X$ has the approximation property, or AP, if the identity
operator on $X$ can be approximated uniformly on compact subsets of $X$
by linear operators of finite rank.  In the 50's, Grothendieck [G]
investigated this property and found several equivalent statements.
For example, he proved that
$X$ has the AP iff the natural map $J:X^*\hat\otimes X\to
X^*\check\otimes X$ is one-to-one ($\hat\otimes$ is the projective
tensor product of Banach spaces and $\check\otimes$ is the injective
tensor product of Banach spaces).  However, it remained unknown if
every Banach space had the AP until Enflo [E] constructed 
the first counter example in the early 70's.
In [S], Szankowski gave a very explicit example of
a subspace of $\ell_p$, $1<p<2$, without the AP.  He considered
$X=(\sum_{n=1}^\infty \oplus \ell_2^n)_p$, which is isomorphic to $\ell_p$, and
defined $Z$ to be the closed span of some vectors of length six.
He then used a clever combinatorial argument to exploit the difference
between the $2$-norm of the blocks
and the $p$-norm of the sum to prove that $Z$ fais the approximation property.  
Szankowski's technique is fairly general. 
In this paper we will use it to show that 
the $\ell_2$-sum (as defined in [P2]) of row operator spaces has a 
subspace without the operator space version of the approximation property.

An operator space $E$ is a Banach space $E$ with an isometric embedding
into $B(H)$, the set of all bounded operators on a Hilbert space $H$.  Or,
equivalently, an operator space $E$ is a closed subspace of
$B(H)$.  If
$E\subset B(H_1)$ and $F\subset B(H_2)$ are operator spaces, their
{\sl minimal tensor product} $E\otimes_{\min} F$ is the closure of the
algebraic tensor product $E\otimes F$ in $B(H_1\otimes_2 H_2)$.  A linear
map $u:E\to F$ is {\sl completely bounded}, or cb, if
for every operator space $G$, the map 
$1_G\otimes u: G\otimes_{\min} E\to G\otimes_{\min} F$ 
is bounded.  The completely bounded norm of $u$, or $\|u\|_{cb}$, is the
supremum of $\|1_G\otimes u\|$, where $G$ runs over all operator spaces
$G$.  It turns out that it is enough to verify that $1_G\otimes u$ is bounded
when $G$ is
$K(\ell_2)$, the set of all compact operators 
on the Hilbert space $\ell_2$, and that 
$\|u\|_{cb}=\|1_{K(\ell_2)}\otimes u\|$.  The set of all cb-maps
from $E$ to $F$ is denoted by $CB(E,F)$. Independently of each other, 
Blecher and Paulsen [BP]
and Effros and Ruan [ER1] gave $E^*$, the Banach space dual of $E$, 
an operator space structure that gives $E^*\otimes_{\min}F$ the norm
induced by $CB(E,F)$.  This indicates that the minimal tensor product is the
operator space analogue of the injective tensor product of Banach spaces.  
In the same papers, Blecher and Paulsen [BP] and Effros and Ruan [ER1]
introduced the operator space analogue of the projective
tensor product.  This is denoted by $E\hat\otimes F$ and satisfies
$(E\hat\otimes F)^*=CB(E,F^*)$.  We refer to [ER3], [J], and [P3]
for more information about operator spaces.

It is well known that the compact subsets of a Banach space $X$ are
contained in the convex hull of null sequences in $X$.  Since there
is a correspondence between null sequences in $X$ and elements of
$c_0\check\otimes X$, it is easy to see that $X$ has the AP iff for every
$u\in c_0\check{\otimes} X$ and every $\epsilon>0$, there exists
a finite rank operator on $X$ such that $\|u-(I\otimes T)(u)\|<\epsilon$.
Based on this observation, Effros and Ruan [EF2] said that an operator
space $V$ has the operator approximation property, or OAP, if for every
$u\in K(H)\otimes_{min} V$ and every $\epsilon>0$, there exists
a finite rank operator $T$ on $V$ such that $\|u-(I\otimes T)(u)\|<\epsilon$.
They proved
that an operator space $V$ has the OAP if and only if the natural map
$J:V^*\hat\otimes V\to V^*\otimes_{min} V$ is one-to-one.  

The following criterion allows us to check that $J$ is not one-to-one,
when $V$ fails the OAP.

\proclaim {Enflo's Criterion} If there exist a sequence of finite rank operators
$\beta_n\in V^*\otimes V$ satisfying:
\item{(i)} $trace(\beta_n)=1$ for every $n\in \N$,
\item{(ii)} $\|\beta_n\|_{V^*\otimes_{min}V}\to 0$ as $n\to\infty$, and
\item{(iii)}$\sum_{n=1}^\infty \|\beta_n-\beta_{n-1}\|_{V^*\hat\otimes V}<\infty$,

\noindent then $V$ does not have the OAP.
\endproclaim

Indeed, $\beta=\beta_1+\sum_{n=2}^\infty \beta_n-\beta_{n-1}=\lim_n\beta_n$
belongs to $V^*\hat\otimes V$ by (iii).  $J\beta=0$ by (ii).  And since
$tr(\beta)=1$, $\beta$ is not zero.  Hence $J$ is not one-to-one and
$V$ fails the OAP.

\head The Construction\endhead

For each $n\in\N$, let
$\Delta_n$ be a partition of 
$\sigma_n=\{2^n, 2^n+1, 2^n+2,\dots, 2^{n+1}-1\}$.
Then
$\{B\in\Delta_n:n\in\N\}$ is a partition of $\N$.  For each
$B\in\Delta_n$, let $R_B$ be the row Hilbert space with orthonormal basis
$\{e_j:j\in B\}$.  
We define $X$ to be the $\ell_2$-sum of these row spaces.  More precisely, 
$X$ is the complex interpolation space between 
$(\sum_{n=1}^\infty\sum_{B\in\Delta_n}\oplus R_B)_\infty$ and 
$(\sum_{n=1}^\infty\sum_{B\in\Delta_n}\oplus R_B)_1$ of parameter 
$\theta={1\over2}$ (see [P2], page 34).  That is, 
$$X=\biggl(\sum_{n=1}^\infty\sum_{B\in\Delta_n}\oplus R_B\biggr)_{\ell_2}=
\left(
\biggl(\sum_{n=1}^\infty\sum_{B\in\Delta_n}\oplus R_B\biggr)_\infty,
\biggl(\sum_{n=1}^\infty\sum_{B\in\Delta_n}\oplus R_B\biggr)_1
\right)_{1\over2}.$$
At the Banach space level, $X$ is a Hilbert space with orthonormal
basis $\{e_i:i\in\N\}$.  But at the operator space level, $X$ is a combination
of row Hilbert spaces and $OH$, the self dual operator Hilbert space
introduced by Pisier in [P1].  If $A\subset\N$, let 
$X_A=\overline{{\text span}}\{e_i:i\in A\}\subset X$. It follows
from the definition of $X$ that if there exists $n\in\N$ such that 
$A\subset B$ for some $B\in\Delta_n$, then 
$X_A$ is completely isometric to $R_A$, the row Hilbert space with 
orthonormal basis $\{e_i:i\in A\}$.  And if for each $n\in\N$,
$A$ has at most one point
from each element in $\Delta_n$,  
(i.e, $\text{card}(A\cap B)\leq 1$ for every
$B\in\Delta_n$, $n\in\N$), then $X_A$ is completely
isometric to $OH_A$.  

\medbreak

Let $Z$ be the closed subspace of $X$ spanned by
$$z_i=e_{2i}-e_{2i+1}+e_{4i}+e_{4i+1}+e_{4i+2}+e_{4i+3},\quad i=1,2\cdots.$$

\proclaim{Theorem 1} With the appropriate selection of $\Delta_n$, $Z$
does not have the OAP.
\endproclaim

For each $i\in\N$, let $z_i^*={1\over 2}(e_{2i}^*-e_{2i+1}^*)$, where
the $e_i^*$'s are biorthogonal to the $e_i$'s.  Then let
$$\beta_n={1\over 2^n}\sum_{i\in\sigma_n} z_i^*\otimes z_i\quad\quad
\text{for }n\geq2.$$ 

We need to check that the $\beta_n$'s satisfy the conditions of Enflo's
criterion.

\bigbreak
\noindent $\underline{ \text{Condition (i)}}$
\smallbreak

This is trivially verified.  Since $z_i^*(z_i)=1$ for every $i\geq1$, we see that
$trace(\beta_n)=(1/2^n)\sum_{i\in\sigma_n}z^*_i(z_i)=(1 / 2^n)|\sigma_n|=1$.

\bigbreak
\noindent $\underline{\text{Condition (ii)}}$
\smallbreak

Since  
$\|\beta_n\|_{Z^*\otimes_{min}Z}\leq \|\beta_n\|_{X^*\otimes_{min}X}
=\|\beta_n\|_{cb}$, we will estimate the cb-norm of
$\beta_n:X\to X$.  However, it follows from the definition of $\beta_n$
that we only need to estimate
the cb-norm of $\beta_n:X_{\sigma_{n+1}}\to X_{\sigma_{n+1}\cup\sigma_{n+2}}$,
where $X_{\sigma_k}=\text{span}\{e_i:i\in\sigma_k\}$.  
Let $I_1:X_{\sigma_{n+1}}\to R_{\sigma_{n+1}}$ and
$I_2: R_{\sigma_{n+1}\cup\sigma_{n+2}}\to X_{\sigma_{n+1}\cup\sigma_{n+2}}$
be the formal identity maps, and let 
$\tilde{\beta}_n:R_{\sigma_{n+1}}\to R_{\sigma_{n+1}\cup \sigma_{n+2}}$
be $\tilde{\beta}_n={1\over 2^n}\sum_{i\in\sigma_n} z_i^*\otimes z_i$
(that is, $\tilde{\beta}_n$ has the same matrix representation of $\beta_n$,
but it is defined on row operator spaces).  Then $\beta_n=I_2\circ
\tilde{\beta}_n\circ I_1$.  Since the $z_i$'s, $i\in\sigma_n$, have
disjoint support, the $z_i^*$'s, $i\in\sigma_n$, have also disjoint support,
and the row spaces are homogeneous, it is easy to see that 
$\|\tilde{\beta}_n\|_{cb}=\|\tilde{\beta}_n\|={1\over 2^{n+1}}\sqrt{3}$.
We will prove condition (ii) by
checking that $\|I_1\|_{cb}\|I_2\|_{cb}\leq\sqrt{2^{n+2}}$.

From the definition of $X$, we see that 
$X_{\sigma_{n+1}}$ is equal to 
$(\sum_{B\in\Delta_{n+1}}\oplus R_B)_{\ell_2}$, the complex
interpolation space $\bigl( (\sum_{B\in\Delta_{n+1}}\oplus R_B)_\infty, 
(\sum_{B\in\Delta_{n+1}}\oplus R_B)_1\bigr)_{1\over 2}$.
It is easy to check that $\|I_1:(\sum_{B\in\Delta_{n+1}}\oplus R_B)_\infty
\to R_{\sigma_{n+1}}\|_{cb}\leq\sqrt{|\Delta_{n+1}}$ and 
that $\|I_1:(\sum_{B\in\Delta_{n+1}}\oplus R_B)_1
\to R_{\sigma_{n+1}}\|_{cb}\leq 1$.  Therefore 
$\|I_1\|_{cb}\leq (|\Delta_{n+1}|)^{1\over 4}$.  Similarly,
$\|I_2\|_{cb}\leq (|\Delta_{n+1}|+|\Delta_{n+2}|)^{1\over 4}$.  Since
$|\Delta_k|\leq 2^k$ for every $k\in\N$, we see that 
$\|I_1\|_{cb}\|I_2\|_{cb}\leq\sqrt{2^{n+2}}$.

\bigbreak
\noindent $\underline{\text{Condition (iii)}}$
\smallbreak

\bigbreak

Using the fact that 
$z^*_i={1\over 4}(e_{4i}^*+e_{4i+1}^*+e_{4i+2}^*+e_{4i+3}^*)$ on $Z$,
we get that
$$
\eqalign{ 
&\beta_n-\beta_{n-1}=
{1\over 2^{n+1}}\sum_{i\in\sigma_n} (e^*_{2i}-e^*_{2i+1})\otimes
			(e_{2i}-e_{2i+1}+e_{4i}+e_{4i+1}+e_{4i+2}+e_{4i+3})\cr
&-
{1\over 2^{n+1}}\sum_{i\in\sigma_{n-1}} 
	(e_{4i}^*\!+\!e_{4i+1}^*\!+\!e_{4i+2}^*\!+\!e_{4i+3}^*)
	\otimes (e_{2i}\!-\!e_{2i+1}\!+\!e_{4i}\!+\!e_{4i+1}
	\!+\!e_{4i+2}\!+\!e_{4i+3})\cr}$$
$$=\scriptstyle 
{1\over 2^{n+1}}\!\!\!\!
\displaystyle \sum_{i\in\sigma_{n-1}}\!\!\!
 \left\{ 
   \eqalign{ 
   & \scriptstyle e^*_{4i}\otimes
   (e_{4i}\!-e_{4i+1}\!+e_{8i}\!+e_{8i+1}\!+e_{8i+2}\!+e_{8i+3}
   \!-e_{2i}\!+e_{2i+1}\!-e_{4i}\!-e_{4i+1}\!-e_{4i+2}\!-e_{4i+3})
   \cr
   & \scriptstyle e^*_{4i+1}\otimes
   (-e_{4i}\!+e_{4i+1}\!-e_{8i}\!-e_{8i+1}\!-e_{8i+2}\!-e_{8i+3}
   \!-e_{2i}\!+e_{2i+1}\!-e_{4i}\!-e_{4i+1}\!-e_{4i+2}\!-e_{4i+3})
   \cr
   &\scriptstyle e^*_{4i+2}\otimes
   (e_{4i\!+\!2}\!-e_{4i\!+\!3}\!+e_{8i\!+\!4}\!+e_{8i\!+\!5}\!+e_{8i\!+\!6}\!+
   e_{8i\!+\!7} \!-e_{2i}\!+e_{2i\!+1}\!-e_{4i}\!-e_{4i\!+1}\!-
   e_{4i\!+2}\!-e_{4i\!+3})
   \cr
   & \scriptstyle e^*_{4i\!+\!3}\otimes
   (-e_{4i\!+\!2}\!+e_{4i\!+\!3}\!-e_{8i\!+\!4}\!-
   e_{8i\!+\!5}\!-e_{8i\!+\!6}\!-e_{8i\!+\!7}
   \!-e_{2i}\!+e_{2i\!+\!1}\!-e_{4i}\!-e_{4i\!+\!1}\!-e_{4i\!+\!2}\!-e_{4i\!+3})
   \cr}
  \right\} 
$$
Note that after cancellation, each of the vectors in the 
parenthesis has nine terms.  Two of them cancel out and two are
equal.  Then we can write each of them as a linear 
combination of nine vector basis.  Eight of them 
have coefficients equal to $\pm1$ 
and the other has a coefficient equal to $\pm2$.

Szankowski defined nine functions $f_k:\N\to\N, k\leq 9$, to index these
vectors.  Let $n=4i+l$ and $l=0, 1, 2, 3$.
Then $f_1(4i+l)=2i$ and $f_2(4i+l)=2i+1$.
For $k=3, 4, 5$, $f_k(4i+l)=4i+[(l+1)\mod 4]$. For $l=0,1$, 
$f_6(4i+l)=8i$, $f_7(4i+l)=8i+1$, $f_8(4i+l)=8i+2$, and $f_9(4k+l)=8i+3$.
And finally, for $l=2, 3$, 
$f_6(4i+l)=8i+4$, $f_7(4i+l)=8i+5$, $f_8(4i+l)=8i+6$, and $f_9(4k+l)=8i+7$.
Then we have
$$\beta_n-\beta_{n-1}=
	{1\over 2^{n+1}}\sum_{j\in\sigma_{n+1}} e_j^*\otimes y_j,$$
where $y_j=\sum_{k=1}^9=\lambda_{j,k}e_{f_k(j)}\in Z$.  Recall that 
eight of the
$\lambda_{j,k}$'s have absolute value equal to one, and one has absolute
value equal to 2.

\bigbreak

The following Lemma of Szankowski provides the key combinatorial
argument (see [S] and [LT, page 108]). 

\proclaim{Lemma 2 (Szankowski)} There exist partitions $\Delta_n$ and
$\nabla_n$  of $\sigma_n$ into disjoint sets, and a sequence 
$m_n\geq 2^{{n\over8}-2}$, $n=2, 3, \dots$, so that 
\item{(1)} $\forall A\in\nabla_n, m_n\leq \text{card}(A)\leq 2m_n$,
\item{(2)} $\forall A\in\nabla_n, \forall B\in\Delta_n, 
			\text{card}(A\cap B)\leq 1$,
\item{(3)} $\forall A\in\nabla_n, \forall 1\leq k\leq 9, f_k(A)$ is contained
in an element of $\Delta_{n-1}, \Delta_n$, or $\Delta_{n+1}$.
\endproclaim

(Notice that $f_k(\sigma_n)\subset\sigma_{n-1}$ for $k=1,2$, 
$f_k(\sigma_n)\subset\sigma_{n}$ for $k=3, 4, 5$, and 
$f_k(\sigma_n)\subset\sigma_{n+1}$ for $k=6, 7, 8, 9$).

Since $\nabla_{n+1}$ is a partition of $\sigma_{n+1}$, we have that
$$\beta_n-\beta_{n-1}={1\over 2^{n+1}}\sum_{A\in\nabla_{n+1}}
\left[\sum_{j\in A} e_j^*\otimes y_j\right].\leqno{(1)}$$

\proclaim{Lemma 3}  For every $A\in\nabla_{n+1}$, 
$\|\sum_{j\in A} e_j^*\otimes y_j\|_{Z^*\hat\otimes Z}
\leq 18\bigl(\text{card}(A)\bigr)^{3\over4}$.
\endproclaim

The last condition of Enflo's criterion follows 
immediatedly from (1), Lemma 2, and Lemma 3.  Indeed,  
$$\eqalign{
\|\beta_n-\beta_{n-1}\|_{Z^*\hat\otimes Z} & \leq
{1\over 2^{n+1}}\text{card}(\nabla_{n+1}) 18\max_{A\in\nabla_{n+1}}
		\text{card}(A)^{3\over4} \cr
&\leq {1\over 2^{n+1}} {2^{n+1}\over m_{n+1}} 18 (2 m_{n+1})^{3\over4}
\leq {36\over m_{n+1}^{1\over4}},
\cr}$$
which is clearly summable.

\medbreak

We only need to prove Lemma 3.  For this, we need 
the result of Pisier (see remark 2.11
of [P1]) that $CB(R_n,OH_n)=S_4^n$, where
$S_4^n$ is the Schatten 4-class.  Consequently, if $S:OH_n\to R_n$, then
$\|S\|_{OH_n\hat\otimes R_n}=\|S\|_{S_{4/3}^n}$.  
In particular, if $I:OH_n\to R_n$
is the formal identity, $\|I\|_{OH_N\hat\otimes R_n}=n^{3/4}$.

\bigbreak

\noindent {\bf Proof of Lemma 3.}  The element
$\gamma=\sum_{j\in A} e_j^*\otimes y_j\in X^*\hat\otimes Z$ 
induces a finite rank map $\gamma:X\to Z$.  The restriction of $\gamma$
to $Z$ is the map $\alpha=\gamma_{|Z}:Z\to Z$, which clearly
satisfies $\alpha=\sum_{j\in A} q(e_j^*)
\otimes y_j\in Z^*\hat\otimes Z$, where $q=(\iota_Z)^*:X^*\to Z^*$ is the
adjoint of the inclusion $\iota_Z:Z\to X$.  Since 
$(Z^*\hat\otimes Z)^*=CB(Z^*,Z^*)$ we have that
$\|\alpha\|_{Z^*\hat\otimes Z}=\sup\{|\langle T, \alpha\rangle|: 
T:Z^*\to Z^*, \quad \|T\|_{cb}\leq1\}$, where $\langle \cdot,\cdot \rangle$
is the trace duality.  

We will see that we can factor $\alpha$
through the formal identity map $I:OH_A\to R_A$, where $R_A$ is 
the row Hilbert space with basis $\{\delta_j:j\in A\}$.
Recall that the projective tensor norm of $I:OH_A\to R_A$, 
viewed as an element of
$OH_A\hat\otimes R_A$,  is
equal to $\bigl(\text{card}(A)\bigr)^{3\over4}$.

Let $\Psi:R_A\to Z$ be the map defined by $\Psi(\delta_j)=y_j$.  We claim
that $\|\Psi\|_{cb}\leq 18$.  Indeed, if $a_j\in B(\H)$ for $j\in A$, 
$$\sum_{j\in A}a_i\otimes y_j=\sum_{j\in A} 
a_j\otimes \sum_{k=1}^9\lambda_{j,k}e_{f_k(j)}=\sum_{k=1}^9\left[\sum_{j\in A} 
\lambda_{j,k}a_j\otimes e_{f_k(j)}\right].$$
It follows from (3) of Lemma 2 that 
$\{f_k(j):j\in A\}\subset B$ for some $B$ in $\Delta_n, \Delta_{n+1}$, or
$\Delta_{n+2}$.  Then the definition of $X$ implies that
the span 
of the $e_{f_k(j)}$'s for $j\in A$ is a row operator
space.  Hence,
$$\biggl\| \sum_{j\in A}a_i\otimes y_j\biggr\|\leq (9)(2)
\biggl\|\sum_{j\in A}a_ja^*_j\biggr\|^{1\over2}
=18\biggl\|\sum_{j\in A}a_j\otimes \delta_j\biggr\|.$$

Let $X_A=\text{span}\{e_j:j\in A\}$.  By (2) of Lemma 2, all the elements
of $A$ belong to different elements of the partition $\Delta_{n+1}$.  This
implies that $X_A$ is completely isometric to $OH_A$.  Let $P_A:X\to X_A$
be the completely contractive projection onto $X_A$, and let 
$I:X_A\to R_A$ be the formal identity.  Then we have that 
$$\alpha=\Psi\circ I\circ P_A \circ \iota_Z.$$

If $T:Z^*\to Z^*$ is completely bounded,
$$\eqalign{
|\langle T, \alpha\rangle|& =|tr(T^*\circ\alpha)|
			=|tr(T^*\circ \Psi\circ I\circ P_A \circ \iota_Z)|
 =|tr(P_A \circ \iota_Z\circ T^*\circ \Psi\circ I)|\cr
&\leq \|P_A \circ \iota_Z\circ T^*\circ \Psi\|_{cb}
	\|I\|_{OH_A\hat\otimes R_A}\leq 18\|T\|_{cb}
	\bigl(\text{card}(A)\bigr)^{3\over 4}.\cr}$$
This finishes the proof of Lemma 3.

\enddocument